The Origin of the Mystical Number Seven in Mesopotamian Culture:

Division by Seven in the Sexagesimal Number System

Kazuo MUROI

§1. Introduction

In Mesopotamian literary works, including hymns, myths, and incantations, the number 7 often occurs in mysterious circumstances where something of religious importance may be indicated. Although we do not yet completely understand the connotations of the mystical number 7 in these literary works, we know that the Sumerians of the third millennium BCE believed the number 7 – out of the many natural numbers – to be special. It seems that various Sumerian words containing the number 7 had taken firm root in their culture by the twenty-second century BCE at the latest. However, as yet we have had no convincing explanation for the Sumerians' adherence to the number 7. Why did they select the number 7 as a mystical or sacred number?

In this paper I shall attempt to clarify the origin of the mystical number 7[1] by examining literary works containing the number 7, as well as some mathematical problems relating to division by 7.

§2. The oldest examples of the mystical number 7

One of the oldest examples of the mystical number 7 is in the *Early Dynastic Proverbs, Collection One* written in the 26th century BCE[2]:

   lul-7 lu-lu "Seven lies are too numerous."

In addition to the lul-7 "seven lies" of this wordplay, we have several instances of a noun plus 7 in certain unusual orthographic texts, which are called UD-GAL-NUN texts, as well as normal orthographic texts:[3]

    UD-UD-7 "seven gods"

    urì-gal-7 "seven (divine) standards"

    ama-ugu$_4$-7 "mother who bore seven (children)"

    šeš-sila$_4$-7 "seven lamb brothers".

Moreover, there is a famous mathematical problem from those days in which division by 7 is performed:[4]

| | |
|---|---|
| 1 še guru$_7$ | 1 guru$_7$ (= 5,20,0,0 sìla) of barley. |
| sìla-7 | 7 sìla (of barley each) |
| 1 lú šu ba-ti | one man received. |
| lú-bi | (The number of) the men are |
| 45,42,51 | 45,42,51 (= 164571). |
| še sìla su-su 3 | 3 sìla of barley is repaid. |

It is evident that the Sumerian scribe of this tablet had carried out the following calculation:

    $5,20,0,0 = 7 \times 45,42,51 + 3$.

Note that the number 7 is especially chosen as a divisor in this exercise because, in the sexagesimal (base 60) number system, division by 7 is more difficult than division by, for example, 2, 3, or 5. The reciprocals of the numbers 2, 3 and 5 are 0;30, 0;20, and 0;12 respectively, but that of 7 is a recurring sexagesimal fraction 0;8,34,17,8,34,17,….

The Sumerian scribes must have noticed that the difficulty of division by 7 is based on this fact, for we know that the numbers 2, 3, and 5 were occasionally called a-rá-gub-ba "normal factor"[5] and 1/7 = 0;8,34,17,… was probably obtained in the Old Babylonian (about 2000–1600 BCE) mathematical texts. In other words, the number 7 is the first natural number whose reciprocal is difficult to grasp in the sexagesimal notation. I believe that this fact led the Sumerians to consider the number 7 to be special or mystical. In the Old Babylonian literary work *Gilgamesh and Huwawa* (Version B),[6] we find the idea that something special will happen at the seventh trial, which seems to support my assumption:

> hur-sag-1-e in-dè-in-bal ᵍⁱˢerin šà-ga-ni nu-mu-ni-in-pàd / hur-sag-2-kam-ma / hur-sag-3-kam-ma / hur-sag-4-kam-ma / hur-sag-5-kam-ma / hur-sag-6-kam-ma / hur-sag-7-kam-ma bal-e-da-ni ᵍⁱˢerin šà-ga-ni im-ma-ni-in-pàd
>
> "He (Gilgamesh) crossed the first mountain with him (Enkidu), but the cedars were not revealed to his heart. The second mountain, the third mountain, the fourth mountain, the fifth mountain, (and) the sixth mountain. When he was crossing the seventh mountain, the cedars were revealed to his heart."

Imitating the above sentences, we may reconstruct the moment the Sumerians had perceived the peculiarity of the reciprocal of 7:

> I calculated the reciprocal of the number 1 and had no problems. The number 2, … , the number 6. When I was calculating the reciprocal of the number 7, a problem (of a recurring fraction) occurred.

§3. Later examples of the mystical number 7

In Gudea Cylinders A and B,[7] written in twenty-second century BCE, we have six examples of the mystical number 7, which clearly show that Mesopotamian mysticism surrounding the number 7 had been established by then.

(1) é-a sá-7 nam-mi-sig$_{10}$ "(After the first, the second, ... , and the sixth,) he placed the seventh square on the temple." (A. 21. 11)

(2) u$_4$-7-kam-ma-ka é-e im-mi-dab$_6$ "By the seventh day he had set (the seven stelae) up around the temple." (A. 23. 4)

(3) na-7 é-e dab$_6$-ba-bi níg lugal-bi-da šà kúš-kúš-dam "The seven stones surrounding the temple are those which take counsel with its lord." (A. 29. 1–2)

(4) dumu-maš-7 $^d$ba-ú-me "(They are) the goddess Bau's seven children, kids." (B. 11.11)

(5) šíta-sag-7 "the seven-headed mace" (B. 7 .12 and13. 21)

(6) u$_4$-7-ne-éš gemé nin-a-ni mu-da-sá-àm árad-dè lugal-e zag mu-da-gub-àm "(In celebration of the temple,) for seven days, a slave woman was equal to her mistress (in status), and a slave was allowed to stand by his master." (B. 17. 19–21)

In the last example, we may find the remote origin of a week of seven days. Moreover, the seventh day of the month was considered to be special as the following Old Babylonian lexical list[8] shows:

u$_4$-é-7 = *se-b*[*u-tum*] "festival of the seventh day"

  $u_4$-é-15 = *ša-pa-*[*at-tum*] "festival of the fifteenth day"

  $u_4$-é-20 = *eš-ru-*[*ú*] " festival of the twentieth day"

  $u_4$-é-25 = *mi-it-*[*hu-rum*] " festival of the twenty-fifth day".

Thus the Babylonians followed the Sumerians' mysticism of the number 7. There are many examples of the mystical number 7 in the literary works of the Old Babylonian period and the later periods, some of which are listed below[9][10].

(7) me-7-bi "the seven divine powers"

(8) abul-kur-ra-7-bi "the seven gates of the underworld"

(9) di-kud-7-bi "the seven judges"

(10) muš-sag-7 "the seven-headed snake"

(11) muš-mah-ka-7 "the seven-mouthed serpent"

(12) $u_4$-7-àm $gi_6$-7-àm "seven days and seven nights"

(13) lú dumu-ni 7-àm "a man whose sons are seven"

(14) ní-te-a-ni 7-kam-ma "his seventh terror"

(15) me-lem$_4$-a-ni 7-kam-ma "his seventh aura"

(16) ne-mur-7-ta "with seven pieces of charcoal"

(17) izi-gar-mè-7 "the seven torches of battle"

(18) íd-ka-7 "the river of the seven mouths"

(19) $u_4$-7 "the seven storms"

(20) mu-7-àm "(for) seven years"

(21) balag-7 "seven harps"

(22) géštug-7 "seven wisdoms"

(23) nimgir-ra 7-na-ne-ne "the seven heralds"

(24) ur-7 "seven dogs"

(25) nimgir-si-7 "the seven paranymphs"

(26) tum$_9$-7-na "the seven winds"

(27) máš-anše-máš-anše-bi 7-me-eš "The herds of cattle are seven."

(28) šagan-7 "the seven flasks"

(29) gi-di-7 "the seven reed flutes"

(30) barag-gal-7 "the seven great daises"

(31) é-7-kam-ma "the seventh house"

(32) nam-tag-bi 7 a-rá 7 nam-tag-bi du$_8$-a "Its sins are 7 times 7, loosen its sins."

(33) muš-eme-7-bi "the snake with seven tongues"

(34) ka-kéš 7 a-rá-2-àm u-me-ni-kéš "Tie twice seven knots."

(35) 7-àm dingir 7-àm-meš 7-àm dingir-hul-a-meš "the seven gods of the world, and the seven evil gods"

(36) 7-àm dingir an-dagal-la-meš "the seven gods of the wide heavens"

(37) 7 abgal íd-da mú-mú-da "the seven sages created in the river"

    Cf. Seven sages of Greece and Seven sages of the bamboo grove.

(38) 7 gi-urì-gal "seven (divine) standards"

(39) níg-na 7-na "seven censers"

(40) $^d$udu-til-meš 7-*šú-nu* "the seven planets (sun, moon, and five planets)"

(41) mul-mul dingir-7-bi "the Pleiades, its seven gods"[11].

§4. The number 7 in mathematical texts

Around 2000 BCE the Semitic Akkadians took over as rulers of Mesopotamia, and they widely absorbed Sumerian culture, including the sexagesimal place value notation which had been invented by the twenty-second century BCE. Therefore, we have many examples of Old Babylonian mathematical texts which show the fact that the Babylonian scribes also adhered to the number 7 and the reciprocal of it. The most striking example is the cuneiform tablet MS 3956 published by J. Fribery in 2007[12], on the obverse and reverse of which four sexagesimal numbers are written down:

25,57,30 (= $2 \times 3 \times 5^2 \times 7 \times 89$)

20,10,25 (= $5^3 \times 7 \times 83$)

3,4,5,4 (= $2^4 \times 7 \times 61 \times 97$)

2,44,3,45 (= $3^3 \times 5^5 \times 7$)

The number 7, the greatest common divisor of the four numbers, is inscribed on the left-hand side of the tablet.

As another example, I cite here YBC 4652, no. 8[13], which is concisely written in Sumerian and presents a simple linear equation involving the number 7:

na$_4$ ì-pàd ki-lá nu-na-tag igi-7-gál ba-zi igi-13-gál ba-zi-ma / ì-lá 1 ma-na sag-na$_4$

en-nam sag-na$_4$ 1 ma-na 15 + 5/6 gín

"I showed a stone (to him, but) I did not weigh it for him. One-seventh was subtracted, one-thirteenth was subtracted, and I weighed (the remainder). It was 1 ma-na (≈ 500g). What was the original (weight of the) stone? The original (weight of the) stone is 1 ma-na 15 + 5/6 gín."

If we denote the original weight of the stone by $x$, the first step of the process – subtracting one-seventh of the original weight – can be written thus:

$$x - \frac{x}{7}$$

The next step – subtracting one-thirteenth of the new weight – can be written thus:

$$\left(x - \frac{x}{7}\right) - \frac{1}{13}\left(x - \frac{x}{7}\right) = 1,0 \text{ (gín)}.$$

To find $x$, the original weight of the stone, we start by multiplying through by 13:

$$13\left(x - \frac{x}{7}\right) - \left(x - \frac{x}{7}\right) = 13,0$$

Simplify to:

$$(13 - 1)\left(x - \frac{x}{7}\right) = 13,0$$

Now 13 − 1 = 12, so divide through by 12 (in sexagesimal 1/12 is 0;5):

$$\left(x - \frac{x}{7}\right) = 13,0 \times 0;5 = 1,5$$

Multiply through by 7 and simplify:

$$(7 - 1)x = 1,5 \times 7 = 7,35$$

As 7 − 1 = 6, divide through by 6 to find $x$ (in sexagesimal 1/6 = 0;10):

$$x = 7,35 \times 0;10 = 1,15;50 \ \text{(gín)}.$$

This solution for the original weight of the stone is not given in the text.

In this way the Babylonian scribes learned how to deal with the reciprocals of seven, 13, and so on, without calculating their sexagesimal values. However, it should be emphasized that there is no mysticism concerning the number 7 in Babylonian mathematics in itself.

§5. Conclusion

In the middle of the third millennium BCE the Sumerians must have noticed that the reciprocal of the number 7, in contrast to the numbers 1, 2, 3, 4, 5, and 6, could not be expressed by a finite sexagesimal fraction but it recurred every three places.

Since the number 7 is the first natural number that has such a property, it was of particular interest to them and became the representative of their number mysticism. The mystical number 7 was, so to speak, born in the rational mind of the Sumerians and grew in the irrational mind of the same people, and it spread all over the world. It is surprising that there are many words containing the number 7 in modern Japanese, for example, 七福神(Seven lucky gods). Thus the mysticism of the number 7, which originated in Sumer, still influences our daily life – mainly in the domain of superstition.

Notes

(1) I was prompted to study the origin of the mystical number 7 by an excellent essay by Prof. Y. Motohashi:

Yoichi Motohashi, Prime Numbers-Your Gems, arXiv: math / 0512143 v1 [math. HO] 7 Dec 2005, esp., p. 8.

(2) Bent Alster, Early Dynastic Proverbs and Other Contribution to the Study of Literary Texts from Abū Ṣalābīkh, *Archiv für Orientforschung* 38 (1993), pp. 1–45. Alster's translation of line 99 in question is "A lie multiplies seven (times)."

(3) For the UD-GAL-NUN texts and the examples cited here see: R. D. Biggs, *Inscriptions from Tell Abū Ṣalābīkh,* 1974, p. 32, Nos. 165, 166, 201, and Figs 2a and 2b in the paper by Alster cited above.

(4) R. Jestin, *Tablettes sumériennes de Šuruppak*, 1937, no. 50. 1 sìla ≈ 1 litre.

(5) K. Muroi, *Studies in Babylonian Mathematics*, No. 3 (2007), pp. 2–8.

(6) The Electronic Text Corpus of Sumerian Literature (= ETCSL): website by Faculty of Oriental Studies, University of Oxford. ETCSL transliteration: c. 1. 8. 1. 5. 1, lines 50–66. The translation is mine.

(7) ETCSL transliteration: c. 2. 1. 7.

(8) M. Civil and others, *Materials for the Sumerian Lexicon* 13 (1971), pp. 257–258. E. Reiner and others, *The Assyrian Dictionary* (= CAD) S (1984), p. 206.

(9) Examples (7)–(31) are taken from ETCSL.

   (7): 1. 4. 1, line 14.

(8): 1. 4. 1, lines 119 and 125.

(9): 1. 4. 1, line 167.

(10): 1. 6. 1, line 63.

(11): 1. 6. 1, line 138.

(12): 1. 7. 4, Segment D, line 3.

(13): 1. 8. 1. 4, line 266.

(14): 1. 8. 1. 5, line 149.

(15): 1. 8. 1. 5, line 199.

(16): 1. 8. 2. 1, line 293 (= 298).

(17): 1. 8. 2. 1, line 471.

(18): 1. 8. 2. 2, line 35.

(19): 1. 8. 2. 2, lines 172 and 189.

(20): 2. 1. 5, lines 92 and 93.

(21): 2. 1. 5, line 200.

(22): 2. 5. 4. 01, line 71.

(23): 4. 06. 1, Segment A, line 49.

(24): 4. 07. 4, line 29.

(25): 4. 07. 4, line 115.

(26): 4. 12. 1, line 17.

(27): 4. 13. 06, line 28.

(28): 4. 15. 2, line 24.

(29): 4. 16.1, line 4.

(30): 4. 16. 1, line 16.

(31): 4. 80. 2, line 127.

For example (32), see:

P. Michalowski, On the Early History of the Ershahunga Prayer, *Journal of Cuneiform Studies*, vol. 39 (1987), pp. 37–48.

(10) Examples (33)–(40) are taken from the bilingual texts cited in CAD, S, pp. 203–204, 230–231.

(11) Example (41): CAD, Z, 1961, p. 50.

The Akkadian equivalents for the Sumerian words are omitted here.

(12) J. Friberg, *A Remarkable Collection of Babylonian Mathematical Texts* (2007), pp. 41–42. Regrettably, Friberg does not notice the importance of this tablet.

(13) O. Neugebauer and A. Sachs, *Mathematical Cuneiform Texts* (1945), pp. 100–103.